\documentclass[12pt]{article}
\usepackage{amsbsy,amssymb,amsmath,amsfonts, epsfig}%, diagrams}

\begin{document}

%%%%%%%
\newtheorem{theorem}{Theorem}[section] % 1st argument is your name for it
\newtheorem{lemma}[theorem]{Lemma}     % 2nd argument is what is printed
\newtheorem{corollary}[theorem]{Corollary}
\newtheorem{proposition}[theorem]{Proposition}
\newtheorem{defn}[theorem]{Definition}

\newtheorem{rem}[theorem]{Remark}

\def \Lip{{\rm Lip\,}}
\def\pf{ \noindent {\bf Proof: \  }}
\newcommand{\qed}{\hfill\vrule height6pt width6pt depth0pt}
\def\endpf{\qed \medskip} \def\colon{{:}\;}

\newcommand{\decomp}{\mathcal{D}}
\newcommand{\decompX}{\{X_i\}_{i=0}^{\infty}}
\newcommand{\decompY}{\{Y_i\}_{i=0}^{\infty}}
\newcommand{\LS}{L}
\newcommand{\RS}{R}
\newcommand{\LSo}{L_{\mathcal{D}_1}}
\newcommand{\RSo}{R_{\mathcal{D}_1}}
\newcommand{\LSt}{L_{\mathcal{D}_2}}
\newcommand{\RSt}{R_{\mathcal{D}_2}}

\newcommand{\cL}{\mathcal{L}}
\newcommand{\cI}{\mathcal{I}}
\newcommand{\cJ}{\mathcal{J}}
\newcommand{\X}{\mathcal{X}}
\newcommand{\Y}{\mathcal{Y}}
\newcommand{\Z}{\mathcal{Z}}
\newcommand{\cW}{\mathcal{W}}
\newcommand{\cS}{\mathcal{S}}
\newcommand{\opX}{\mathcal{L}(\X)}
\newcommand{\opY}{\mathcal{L}(Y)}
\newcommand{\opW}{\mathcal{L}(W)}
\newcommand{\idI}{\mathcal{I}(\mathcal{X})}
\newcommand{\AD}{\mathcal{A}}
\newcommand{\TD}{T_{\mathcal{A}}}
\newcommand{\N}{\mathbb{N}}
\newcommand{\CN}{\mathbb{C}}

\newcommand{\M}{\mathbb{M}}

\newcommand{\LCN}{\lambda\in\mathbb{C}}
\newcommand{\linf}{\ell_{\infty}}
\newcommand{\oplinf}{\mathcal{L}(\ell_{\infty})}
\newcommand{\ssop}[1]{\mathcal{S}(#1)}
\newcommand{\sso}{\mathcal{S}(\ell_{\infty})}
\newcommand{\comp}[1]{\mathcal{K}(#1)}
\newcommand{\Trest}[1] {\lambda_{#1}I_{|#1} + K_{#1}}
\newcommand{\seq}[1] {\displaystyle \{#1_i\}_{i=1}^{\infty}}
\newcommand{\seqo}[1] {\displaystyle \{#1_i\}_{i=0}^{\infty}}
\newcommand{\sumspace}{\big (\sum_{i=0}^{\infty} Y_i\big )_{p}}
\newcommand{\sumY}{\big (\sum Y\big )_p}
\newcommand{\sumYI}{\big (\sum Y\big )_\infty}
\newcommand{\sumX}{\big (\sum \X\big )_p}
\newcommand{\opsumY}{\mathcal{L}\big (\big (\sum Y\big )_p\big )}
\newcommand{\ILT} {I- \lambda T}
\newcommand{\les}{\sigma_{l.e.}(T)}
\newcommand{\IX}{\mathcal{M}_\X}
\newcommand{\propP}{{\textbf{P}}}

%\addtolength{\jot}{14pt}

\title {{Ideals in $L(L_1)$}
}

\author{ W.~B.~Johnson\thanks{Supported in part by NSF DMS-1565826, Clay Mathematics Institute, and MSRI.},
G.~Pisier\thanks{Supported in part by  MSRI.}, and G.~Schechtman\thanks{ Supported in part by  ISF  and MSRI.}}

\maketitle

%%----------------------------------------------------------------------------------------------------------------------------------------------------

\begin{abstract}
\noindent
The main result is that  there are infinitely many; in fact, a continuum; of closed (two-sided) ideals in the Banach algebra $L(L_1)$ of bounded linear operators on $L_1(0,1)$.  This answers a question from A.~Pietsch's 1978 book ``Operator Ideals". The proof also shows that $L(C[0,1])$ contains a continuum of closed  ideals.  Finally,  a duality argument yields that $L(\ell_\infty)$   has a continuum of closed ideals.
\end{abstract}

\section{Introduction}\label{intro}
After $C^*$-algebras, the spaces of bounded linear operators $L(X)$ on non Hilbertian classical Banach spaces $X$  are arguably the most natural non commutative Banach algebras. In 1969, Berkson and Porta wrote a foundational paper \cite{bp} on $L(X)$ with special attention given to $X=L_p$ ($L_p:= L_p(0,1)$). Previously known was that the ideal of  compact operators $K(X)$ in $L(X)$ is the only non trivial closed ideal when $X$ is one of the classical spaces $\ell_p$, $1\le p < \infty$ or $c_0$ \cite{calkin}, \cite{gohbergetal}, \cite{whitley}. (Throughout this paper ``ideal" means ``two-sided ideal".) In the 1970s there was progress in constructing new  closed ideals in $L(L_p)$.  For example, after Rosenthal \cite{ros} constructed a few non obvious ones,  the third author \cite{schecht} proved that $L(L_p)$ contains infinitely many closed ideals when $1< p \not=2 < \infty$.  Several years later,   Bourgain, Rosenthal, and the third author  \cite{brs}  proved that  there are at least $\aleph_1$.  Actually, these last two results are not even stated in the cited papers, but, as pointed out by Pietsch   \cite[Chapter 5]{pietsch}, it is an easy consequence of the  main results of \cite{schecht} and  \cite{brs} that there are at least $\aleph_0$; respectively, $\aleph_1$;  isomorphically distinct infinite dimensional  complemented subspaces of  $L_p$, $1 < p \not=2 < \infty$, each of which is isomorphic to its square.  After the 1970s, relatively little research was done on the   Banach algebra structure of $L(X)$ spaces until the current millenium. Some of that research culminated in a paper of Schlumprecht and Zsak \cite{sz} in which they prove that there are a continuum of closed ideals in $L(L_p)$ when $1 < p \not= 2 < \infty$. This solved a problem from Pietsch's 1978 book \cite[Chapter 5]{pietsch}.  

 It is perhaps surprising that $L(X)$ with $X$ separable can contain  $2^{2^{\aleph_0}}$ closed ideals. It was pointed out on MathOverFlow by Kania \cite{kania}   that  an example due to Mankiewicz \cite{man} has this  property; in fact, there are even $2^{2^{\aleph_0}}$ maximal  ideals in his example. Several months after we wrote the first version of this paper, the first and third authors constructed  $2^{2^{\aleph_0}}$ closed ideals in $L(L_p)$, $1<p\not= 2 < \infty$. These new ideals are not maximal--it was proved in \cite[Section 9]{jmst} that $L(L_p)$ has a unique maximal ideal. 

\noindent
The situation for $L(L_1)$ is different.  The previously known closed ideals in this Banach algebra are the compact operators $K(L_1)$, which is the smallest one because $L_1$ has the approximation property; the strictly singular operators, ${\mathcal S}(L_1)$;  the  operators on $L_1$ that factor through $\ell_1$; the Dunford--Pettis operators--that is, the operators that map weakly compact sets onto norm compact sets;  and the unique maximal ideal, ${\mathcal M}(L_1)$.  Some explanation is necessary. Given an operator $T:X\to Y$ between Banach spaces and a Banach space $Z$, the operator $T$ is   $Z$--{\sl  singular} if $TS$ is not an (into) isomorphism for any operator $S:Z \to X$.  An operator is strictly singular if it is $Z$--singular for every infinite dimensional $Z$. The strictly singular operators, ${\mathcal S}(X)$, on any Banach space $X$ form a closed ideal, as do the weakly compact operators, ${\mathcal W}(X)$.  But ${\mathcal W}(L_1) = {\mathcal S}(L_1)$ because $L_1$ has the Dunford-Pettis property \cite[Theorm 5.4.6]{albiackalton}.  If $Z$ is any Banach space, ${\mathcal I}_Z(X)$ denotes the operators on $X$ that factor through $Z$. Obviously, $L(X)\cdot {\mathcal I}_Z(X) \cdot L(X) \subset L(X)$, so ${\mathcal I}_Z(X)$ is a (usually non closed) ideal in $L(X)$ if it is closed under addition, which it will be if $Z \oplus Z$ is isomorphic to a complemented subspace of $Z$.  It happens that  the ideal ${\mathcal I}_{\ell_1}(L_1) $ is closed in $L(L_1)$.  This is because it is the same as the Radon--Nikodym operators
\cite[Theorem 24.2.7 and Proposition 24.2.12]{pietsch} by a result of Lewis and Stegall \cite{ls}. The ${\mathcal I}_{\ell_1}(L_1) $ is the smallest ideal that is not contained in the ideal ${\mathcal S}(L_1)$ of strictly singular operators because the identity on $\ell_1$ factors through every non strictly singular operator on $\ell_1$ \cite[Section 5.2]{albiackalton}.  ${\mathcal M}(X)$ denotes the set of operators $T$ in $L(X)$ such that the identity on $X$ does not factor through $T$. Evidently ${\mathcal M}(X)$  is an ideal if it is closed under addition, in which case it is obviously the largest ideal in $L(X)$. Then it must be closed because the invertible elements in any Banach algebra form an open set.  Enflo and Starbird \cite{es} proved that ${\mathcal M}(L_1)$ is closed under addition and that ${\mathcal M}(L_1)$ is the same as the $L_1$--singular operators on $L_1$. That the ideal of Dunford--Pettis operators is different from the other four mentioned ideals is due to Coste; the proof is in \cite[p. 93]{du}.
In his  book \cite{pietsch}, Pietsch asked  whether there are infinitely many closed ideals in  $L(L_1)$.  Until now no one has proved that $L(L_1)$ contains a closed ideal different from the four mentioned above. (We are indebted to T.~Kania for sending us unpublished notes that made writing this paragraph easier.)

\noindent
A common way of constructing a (not necessarily closed) ideal in $L(X)$ is to take some operator $U:Y\to Z$ between Banach spaces and let ${\mathcal I}_U(X)$ be the collection of all operators on $X$ that factor through $U$; that is, all $T\in L(X)$ such that there are  $ A\in L(X,Y)$ and $B\in L(Z,X)$ such that  $T = BUA$.  We write ${\mathcal I}_Z(X)$ for ${\mathcal I}_{I_Z}(X)$ as we did above. $L(X){\cal I}_UL(X)\subset {\cal I}_U$ is clear, so ${\cal I}_U$ is a (proper) ideal in $L(X)$ if ${\cal I}_U$ is closed under addition and $I_X$ does not factor through $U$. One usually guarantees the closure under addition by using a $U$ for which  $U\oplus U: Y\oplus Y \to Z\oplus Z$ factors through $U$, and these are the only $U$ that will  appear in the sequel.  Each of the new ideals in $L(L_1)$ that we construct is of the form $\overline{{\mathcal I}_U(L_1)}$ for some operator $U: \ell_1 \to L_1$. Using the standard construction technique we mentioned,  it is easy to build a continuum of ideals in $L(L_1)$; the difficulty is to show that the closure of the ideals are different. To illustrate this difficulty, consider, for example, the family ${\cal I}_{L_p}(L_1)$, $1<p<2$. These are all different, but their
closures $\overline{{\cal I}_{L_p}(L_1)}$ are all equal to  ${\mathcal S}(L_1)$. Since this is not really relevant for our main result, we only outline a proof for those who know the relevant background in the local theory of Banach spaces. To see that these ideals are different, fix $1<p<q<2$. Let $T$ be a surjection from $L_1$ onto the closed span of a sequence of IID symmetric $p$-stable random variables in $L_1$. So the image of $T$ is isometrically isomorphic to $\ell_p$.  Were $T$ to factor through $L_q$, the space $\ell_p$ would be isomorphic to a quotient of a subspace of $L_q$, which it is not because, for example, quotients of subspaces of $L_q$ have type $q$ and $\ell_p$ does not.
Next we show that $\overline{{\cal I}_{L_p}(L_1)}$ is the ideal of weakly compact operators when $1<p<\infty$.  The containment of $\overline{{\cal I}_{L_p}(L_1)}$ in ${\mathcal W}(L_1)$ is clear. Let $T\in {\mathcal W}(L_1)$ and let $\epsilon >0$. By the classification of weakly compact sets in $L_1$ we have that there is a constant $M< \infty$ so that $T_{B_{L_1}} \subset M {B_{L_\infty}} + \epsilon {B_{L_1}} \subset M {B_{L_p}} + \epsilon {B_{L_1}} $. Take an increasing sequence $E_n$ of subspaces of $L_1$ so that $E_n$ is isometrically isomorphic to $\ell_1^n$ and $\cup_{n\in \N} E_n$ is dense in $L_1$. Let $P_n$ be a contractive projection from $L_1$ onto $E_n$. Choose $x_{n,i} $ in $M {B_{L_p}}$ so that $\|T e_{n,i} - x_{n,i}\|_1 \le \epsilon$, where $(e_{n,i} )_{i=1}^n$ in $E_n$ is isometrically equivalent to the unit vector basis of $\ell_1^n$. Let $T_n:L_1 \to L_1$ be $P_n$ followed by the linear extension to $E_n$ of the mapping $e_{n,i} \mapsto x_{n,i}$.  Since  $T_nB_{L_1} $ is a subset of the weakly compact set $ M {B_{L_p}}$, a subnet of $T_n$ converges in the weak operator topology to an operator $S$ on $L_1$. By construction, $\|S x - Tx\|_1 \le \epsilon \|x\|$ if $x\in \cup_{n\in \N} E_n$ and hence by density  for all $x\in L_1$.

\noindent
We now describe without proof a somewhat less elementary example where different ideals have the same closure.
Let $U:\ell_1 \to L_1$ be an injective operator that maps the unit vector basis for $\ell_1$ onto the  Rademacher functions. Several months before completing the research herein, we proved that $\overline{{\mathcal I}_U(L_1)}$ is an ideal different from the four previously known ones.  It was natural then to look at the ideals ${\mathcal I}_{U_p}(L_1)$, $1<p<2$, where $U_p:\ell_1 \to L_1$ is an injective operator that maps the unit vector basis for $\ell_1$ onto an IID sequence of symmetric $p$--stable random variables.  These ideals are all distinct, but it turned out that for all $p$ we have $\overline{{\mathcal I}_{U_p}(L_1)} = \overline{{\mathcal I}_U(L_1)}$.

\noindent
It is convenient to break ideals in $L(X)$ into two classes.  An ideal in $L(X)$ is {\sl small} if it is contained in the strictly singular operators; otherwise it is called {\sl large}. The space $L(L_1)$ is unusual in that every large closed ideal contains the the ideal of strictly singular operators, which in every space is the largest small ideal. In fact, the ideal in $L(L_1)$ of operators that factor through $\ell_1$ is the smallest large ideal and contains the ideal of strictly singular operators.
Obviously  ${\mathcal I}_U(X)$ is small if $U$ is strictly singular, and ${\mathcal I}_U(X)$  is large if $U:Y\to Z$ maps an isomorph of an infinite dimensional complemented subspace of $X$ onto a complemented subspace of $Z$; for example, $U=I_Y$ for some complemented subspace $Y$   of $X$.  The first  condition is not necessary for smallness.       Indeed, if $Z=\ell_q$ for some $ 2<q<\infty$, then every operator from $\ell_q $ into $L_1$ is compact \cite{ros2}, so ${\mathcal I}_{I_Z}(L_1)$ is small although $I_Z$ is not strictly singular.  Similarly, the second condition  is not necessary for largeness. Consider, for example, $X = C[0,1] \oplus \ell_2$. Let $T$ be an isometric embedding of $\ell_2$ into $C[0,1]$ and define $U:X\to X$ by $U(x,y) := (Ty,0)$.  Certainly $U$ is not strictly singular, but for any infinite dimensional subspace $Y$ of $X$, $UY$ is not complemented in $X$. Indeed,  $U$ is weakly compact,
but no infinite dimensional reflexive subspace of $C[0,1]$ is complemented because $C[0,1]$ has the Dunford--Pettis property \cite[Theorem 5.4.6]{albiackalton}.
A more striking example was given by    Astashkin, Hern\'andez,  and Semenov \cite{ahs}. Fix $0<\lambda <1/2$ and let $J$ be the formal identity operator from $L \log^\lambda L(0,1)$ into $L_1(0,1)$. This operator is an isomorphism on the closed span of the Rademacher functions, but Astashkin, Hern\'andez,  and Semenov  \cite{ahs}   proved that $J$ is not an isomorphism on any infinite dimensional complemented subspace of $L \log ^\lambda L$. Let $X:= L_1 \oplus L \log^\lambda L$ and define $U: X \to X$ by $U(x,y)  = (Jy,0)$. Then $U$ is not strictly singular but also is not an isomorphism when restricted to any  infinite dimensional complemented subspace of $X$.

\noindent
It clarifies what we already said about known results to mention that the results in   \cite{brs}   show that there are at least $\aleph_1$ large closed ideals in $L(L_p)$, $1 < p \not=2 < \infty$, while in \cite{sz} it was proved that there are at least a continuum of small closed ideals  in these spaces. The $2^{2^{\aleph_0}}$ closed ideals constructed recently by the first and third authors are large, but it is open whether there are more than a continuum of small closed ideals in $L(L_p)$, $1 < p \not=2 < \infty$.  In Section \ref{ideals} we prove that there are a continuum of small closed ideals  in $L(L_1)$. We do not know whether there are more than three large ideals in $L(L_1)$. This is related to the famous problem whether every infinite dimensional complemented subspace of $L_1$ is isomorphic to $\ell_1$ or $L_1$. 
Also, in contrast to the case of $p>1$, we do not know if there are more than a continuum of closed ideals in $L(L_1)$.  

An immediate consequence of the construction in Section  \ref{ideals} is that $L(C(\Delta))$, where $\Delta := \{-1,1\}^{\N}$ is the Cantor set,  contains a continuum of  small closed ideals.  By Miljutin's Theorem \cite{mil},  \cite[Theorem 4.4.8]{albiackalton}, this is the same as saying that $L(C(K))$ has a continuum of small closed ideals for every compact uncountable metric space $K$.  
In Section \ref{ellinfinity} we use a duality argument to prove that $L(\ell_\infty)$ has a continuum of small closed ideals and also show more generally that distinct small  closed ideals in $L(L_1)$ dualize to give distinct small closed ideals in $L(L_\infty)$.  When $X$ is reflexive, it is obvious that distinct closed ideals in $L(X)$ dualize to give distinct closed ideals in $L(X^*)$, but one does not expect this to be the case when $X$ is non reflexive.

\noindent
Our notation is standard. We just mention that when a quantity (usually of the form $N^\alpha$ with $N\in \N $ but $\alpha \not\in \N$) is treated as an integer, it should be adjusted to the closest larger or smaller integer, depending on context. As we mentioned earlier, all ideals are assumed to be two-sided ideals.  All results are valid for both real and complex Banach spaces even if some proofs are, for simplicity of notation, written for the case of real scalars. 

\noindent
Much of this research was conducted while the authors participated in the Mathematical Sciences Research Institute  program
``Geometric Functional Analysis and Applications" in the Autumn of 2017. 
The authors thank MSRI  for support during that semester.

\section{A Lemma}
\label{thelemma}
In this section we prove a   lemma that will be used in the construction of a continuum of closed ideals in $L(L_1)$. Although the proof is quite simple, it was surprising to us and we think that it will be useful again down the road.  In the lemma, the domain $L_1$ space should be $L_1$ of a probability space (else condition (1) must be adjusted); we wrote the proof for $L_1(0,1)$ with  Lebesgue measure. As for the range space, it follows   formally from the lemma as stated that the range $L_1$ space can be any $L_1$ space of dimension ${N^{\frac{p}{2}}}$--it is just more convenient to prove it for $L_1^{N^{\frac{p}{2}}} $, where the measure is the uniform probability measure on $ {N^{\frac{p}{2}}}$ points.

\begin{lemma}\label{lemma}
Let $1\le p < q<\infty$,  $\{v_1,\dots,v_N\} \subset L_q$, and let $T:L_1 \to L_1^{N^{\frac{p}{2}}}$ be an operator.  Suppose that $C$ and $\epsilon$ satisfy
\begin{enumerate}
\item $ \max_{|\epsilon_i|=   1}\| \sum_{i=1}^N \epsilon_i v_i \|_q \le C N^{1/2}$, and
\item $ \min_{1\le i \le N} \| Tv_i \|_1 \ge \epsilon$.
\end{enumerate}
Then  $\|T\| \ge ( \epsilon/C) N^{{q-p}\over{2q}}$.
\end{lemma}
\pf
Take $u_i^* $ in $L_\infty^{N^{p/2}} = (L_1^{N^{\frac{p}{2}}} )^*$ with $|u_i^*| \equiv 1$ so that $\langle u_i^*, Tv_i \rangle = \|T v_i\|_1 \ge \epsilon$.  Then
%\begin{equation}
\begin{align*}%{rl}
\epsilon N & = \sum_{i=1}^N \langle T^* u_i^*, v_i \rangle :=  \int_0^1  \sum_{i=1}^N (T^* u_i^*)(a) v_i(a) \, da \\
& \le \int_0^1   \sup_{a\in [0,1]} |    \sum_{i=1}^N (T^* u_i^*)(a) v_i(b)| \,db \\
&=: \int_0^1 \|  \sum_{i=1}^N v_i(b) T^* u_i^*\|_{{L_\infty [0,1]}^{\phantom{N^N}}} db \\
&\le \|T\|  \int_0^1 \|  \sum_{i=1}^N v_i(b)  u_i^*\|_{L_\infty^{N^{p/2}}} \,db \\
&\le \|T\| N^{\frac{p}{2q}}  \int_0^1   \big{(}\int_{[N^{\frac{p}{2}}]} |  \sum_{i=1}^N u_i^*(c)  v_i(b) |^q \,dc\big{)}^{\frac{1}{q}} \, db \\
&\le  \|T\| N^{\frac{p}{2q}} \big{(}     \int_{[N^{\frac{p}{2}}]}     \int_0^1   |  \sum_{i=1}^N u_i^*(c)  v_i(b) |^q \,db \, dc\big{)^{\frac{1}{q}}} \\
&\le C \|T\| N^{\frac{p+q}{2q}}. \quad\quad  \quad\quad \quad\quad  \quad\quad \quad\quad \quad\quad \quad\quad \endpf
\end{align*}
%\end{equation}
\begin{rem}  For $q\ge 2$,
the power of $N$ in the conclusion of Lemma \ref{lemma} is of optimal order as $N\to \infty$.
\end{rem}
This is shown by the following simple argument.
Set $m:=N^{p/2}$. By    \cite{bgn}
there are functions $f_1,\cdots,f_{m^{2/q}}$
that are unit vectors in $L^m_1$ (these functions are even $\pm 1$ valued) 
such that  $$\|\sum x_j f_j\|_{L_q^m} \le b_q (\sum|x_j|^2)^{1/2}$$
for every $x=(x_j) \in \ell_2^{m^{2/q}}$.
Here $b_q $ is a constant depending only on $q>2$.
Assume $p\le q$ so that $N \ge {m^{2/q}}$. 
Let  $K$  be (the closest integer to)  $m^{-2/q} N$.  
Let $s: [N] \to [m^{2/q}]$
be a surjection such that $|s^{-1}(k)|= K$
for every  $k\in [m^{2/q}]$.
Define   $v : \ell_\infty^N \to L_q^m $   by
$v(x)= \sum x_k f_{s(k)}$ for all $x\in \ell_\infty^N$.
We claim that
$\|v\|\le   b_q  K {m^{1/q}}$.
Indeed, $$\| v(x)\|_q\le \| \sum (\sum\nolimits_{j\in s^{-1}(k)}  x_j) f_j\|_q
\le   b_q ( \sum |\sum\nolimits_{j\in s^{-1}(k)}  x_j|^2 )^{1/2} \le b_qK ( {m^{2/q}})^{1/2} \sup|x_j|. $$
Let $T: L_1^m \to L_1^m $ be the identity and
let $J_m: L_q^m \to L_1^m$ be the inclusion.
We set $v_k:= J_m v(e_k) \in L_1^m$. Note  that  $\|v_k\|_1=1$
and hence
$\epsilon=\inf \|v_k\|_1=1$.
Then condition 1. in Lemma \ref{lemma}
holds with $C:=b_q  K {m^{1/q}} N^{-1/2}$.
We find
$$1=\|T\|= \epsilon  $$
and also
$$( \epsilon/C) N^{{q-p}\over{2q}} =C^{-1} N^{1/2} m^{-1/q} =  b^{-1}_q Nm^{-2/q}K^{-1} = b^{-1}_q  .$$
Thus if the conclusion of the Lemma
is  $\|T\| \ge ( \epsilon/C) N^{\alpha}$
 with an exponent $\alpha$ we have necessarily $\alpha \le {{q-p}\over{2q}}$.
This proves our claim.
 \endpf

\section{Small Ideals in $L(L_1)$}
\label{ideals}

First we fix some notation that will be used in the proof of Theorem \ref{theorem}. Let $\mu$ be the uniform probability on the two point set $\{-1,1\}$. Given a set $S$, $L_1\{-1,1\}^S$ denotes the space of $\mu^S$-integrable functions on $\{-1,1\}^S$. If $T\subset S$, we regard $L_1\{-1,1\}^T$ as a subspace of $L_1\{-1,1\}^S$, where $f$ in $L_1\{-1,1\}^T$ is identified with $\tilde{f}$ in $L_1\{-1,1\}^S$, defined by $\tilde{f}(x_t)_{t\in S} := f(x_t)_{t\in T}$. Our model for $L_1$ is $L_1\{-1,1\}^{\N}$.

\begin{theorem}\label{theorem}
There exists a family $\{{\cal{I}}_p : 2<p<\infty \}$ of (non-closed) ideals in $L(L_1)$ such that their closures $\overline{{\cal{I}}_p}$ are distinct small ideals in $L(L_1)$.
\end{theorem}
\pf
Write $\N$ as the disjoint union of finite sets $E_k$ so that for each $n \in \N$, the cardinality $|E_k|$ of $E_k$ is $n$ for infinitely many $k$. For $n\in \N$, fix $k_n$ so that $|E_{k_n}| = n$.  Fix $2<p<\infty$, fix $n\in \N$, and define $N:= N(p,n)$ to be $2^{2n/p}$, so that the dimension of $L_1\{-1,1\}^{E_{k_n}} $ is $ N^{p/2}$. In view of Bourgain's \cite{bou} solution to Rudin's $\Lambda (p)$-set problem, for some constant $C_p<\infty$ there are non constant (and hence mean zero)   characters $\{v_i(k_n,p) : 1\le i \le N\}$ on the group  $\{-1,1\}^{E_{k_n}} $ so that for all scalars $(a_i)_{i=1}^N$,
\begin{equation}\label{eq1}
\|\sum_{i=1}^N a_i v_i(k_n,p) \|_p \le C_p (\sum_{i=1}^N |a_i|^2)^{1/2}.
\end{equation}
(One can substitute  an older and simpler  theorem for Bourgain's deep result--see the remark following the proof.)
For $k\not= k_n$ in $\N$ with $|E_k|= n$, let $f_k: E_{k_n} \to E_k$ be a bijection and let $T_k : {\mathbb K}^{E_{k_n}} \to {\mathbb K}^{E_k}$ be the induced linear isomorphism (${\mathbb K}$ is the scalar field, either $\mathbb C$ or $\mathbb R$). So for all $r$, $T_k$ is an isometric isomorphism from $L_r\{-1,1\}^{E_{k_n}} $ onto $L_r\{-1,1\}^{E_k}$.  Next set $v_i(k,p):= T_k v_i(k_n,p)$ for $1\le i \le N$.  Finally, define, $J_p: \ell_1 \to L_1\{-1,1\}^{\N}$ to be  an injective operator that maps the unit vector basis of $\ell_1$ onto the following set of characters:
\[
V_p := \bigcup_{n=1}^\infty \bigcup_{|E_k|=n}  \{ v_i(k,p): 1 \le i \le N(p,n) \}.
\]
We can now define the ideals ${\cal{I}}_p$, $2<p<\infty$. Set
\[
{\cal{I}}_p := \{T\in L(L_1) :  T \ \text{factors through} \  J_p \},
\]
where, as we mentioned earlier, $L_1:= L_1\{-1,1\}^{\N}$.  Obviously $L(L_1) {\cal{I}}_p  L(L_1) \subset { \cal{I}}_p $. Moreover, by construction, ${\cal{I}}_p $ is closed under addition. This follows from the observation that $J_p \oplus J_p : \ell_1 \oplus \ell_1 \to L_1 \oplus L_1 $ factors through $J_p$. We leave the checking of the observation to the reader, just remarking that for ${\mathbb M}\subset \N$, if you identify $L\{-1,1\}^\M$ with the functions in $L_1\{-1,1\}^{\N}$ that depend only on $\M$, then the (norm one) conditional expectation projection $P_\M $ from $L_1\{-1,1\}^{\N}$  onto $L_1\{-1,1\}^{\M}$ is zero on the mean zero functions in $L_1\{-1,1\}^{\N \setminus \M}$.

\noindent
We show that $\overline{{\mathcal I}_p} \not= \overline{{\mathcal I}_q} $  when $p\not= q$ by verifying that ${{\mathcal I}_p} \not\subset \overline{{\mathcal I}_q} $ when $2 < p < q < \infty$.  We will use the observation that
for all $x_1,\dots,x_n$ in $V_p$ and scalars $a_1,\dots, a_n$, we have $\|\sum_{i=1}^n a_i x_i \| \le  2 B_p C_p  ( \sum_{i=1}^n |a_i|^2)^{1/2}$, where
$B_p$ is the constant in Khintchine's inequality. 
% $V_p$  satisfies (\ref{eq1}) with   constant $C_p'\le  2 B_p C_p$, where $B_p$ is the constant in Khintchine's inequality.  
 Indeed, sequences of mean zero independent random varibles are unconditional in $L_p$ with constant at most $2$, so if $f_k $ is in the linear span of $\{v_i(k,p): 1 \le i \le N(p,k )\}$, then
\[
\| \sum_k a_k f_k\|_p \le 2 B_p \|(\sum_k |a_k|^2 \|f_k\|_p^2)^{1/2}
\] because $L_p$ has type $2$ with constant $B_p$, hence the observation follows from (\ref{eq1}).

\noindent
Suppose $2<p<q<\infty$. We prove that ${{\mathcal I}_p} \not\subset  \overline{{\mathcal I}_q} $ by showing that $J_pQ \not\in \overline{{\mathcal I}_q} $, where $Q$ is a quotient mapping from $L_1$ onto $\ell_1$.
%
%Suppose, for contradiction, that ${{\mathcal I}_p} \subset  \overline{{\mathcal I}_q} $ for some $2 < p < q < \infty$.
%
In fact, we prove the formally stronger fact that $V_p$, the image under $J_p$ of the unit vector basis of $\ell_1$, has  distance at least $1/10$  from    $SB_{L_1}$ for every $S$    in $\cI_q$.
Let $\epsilon > 0$; $ \epsilon =1/10$ will do. If our claim is false,  there is $T$ in $L(L_1)$ so that
\begin{equation}\label{eq2}
V_p \subset TJ_q B_{\ell_1} + \epsilon B_{L_1}^o,
\end{equation}
where $B_X$ is the closed unit ball of $X$ and $B_X^0$ is the open unit ball of $X$. For $n$ in $\N$ define $T_n = P_{E_{k_n}} T $, where again $P_{E_{k_n}} $ is the conditional expectation projection from $L_1\{-1,1\}^\N$ onto $L_1\{-1,1\}^{E_{k_n}}$.  Clearly $\|T_n\| \le \|T\|$. Let $V_p(k_n) := \{v_i(k_n,p) : 1 \le i \le N(p, n )\}$. Then
\begin{equation}\label{eq3}
V_p(k_n) \subset T_n J_q B_{\ell_1} + \epsilon B_{L_1^{E_{k_n}}}^o = \text{conv} T_n(\pm V_q ) + \epsilon B_{L_1^{E_{k_n}}}^o.
\end{equation}
Let
\[
V_q':= \{w\in V_q: \|T_n w\| \ge 1-2\epsilon \}.
\]
We need a lower bound on the cardinality of  $V_q'$; namely, that $|V_q'| \ge \delta N(p,n)$ for some constant $\delta$ that does not depend on $n$, for then by Lemma \ref{lemma}   we get a contradiction by letting $n\to  \infty$.

\noindent
To show that $|V_q'| \ge \delta N(p,n)$, we use the fact that the modulus  $|T_n|$ of the $L_1$ operator $T_n$  has the same norm as $T_n$.  Let $f:= |T_n| 1$.  Then $\|f\|_1 \le \|T\|$ and for all $w$ in $V_q$ we have $|T_n w| \le f$. Define
$E:= [f \le \|T\|/\epsilon] \subset \{-1,1\}^{E_{k_n}}$, so  the measure of its complement $\widetilde{E}$ is at most $\epsilon$. Since each
$v_i  := v_i(k_n,p)   $    has constant modulus $1$, we have that
\begin{equation}\label{eq3}
\langle 1_E v_i, v_i\rangle \ge 1-\epsilon \quad \text{for all} \ 1 \le i \le N(p,n).
\end{equation}
Since $\|v_i - u\|_1 < \epsilon$ for some $u$ in $\text{conv} \, \pm T_n V_q$, we get from (\ref{eq3}) that for each
$1 \le i \le N(p,n)$, there is $w_i$ in $ \pm V_q$ such that $\langle 1_E v_i, T_n w_i \rangle \ge 1 - 2 \epsilon$ and hence $w_i$ is in $V_q'$; that is $\|T w_i \| \ge 1 - 2 \epsilon$. That would complete the proof except for the annoying fact that these $w_i$ need not be distinct. However, for any $w$ in $ V_q$, $w_i = \pm w$ for only a bounded number  of $i$'s, where the bound depends only on $\epsilon$ and $\|T\|$.  Indeed, if $   \langle  1_E v_i, \epsilon_i T_n w \rangle \ge 1 - 2\epsilon$ for some $\epsilon_i = \pm1$ and all $i$ in $S$, then
\begin{align}\label{S}
|S| (1-2\epsilon) &\le \langle 1_E \sum_{i\in S} \epsilon_i  v_i, T_n w \rangle =  \langle  \sum_{i\in S} \epsilon_i  v_i, 1_E T_n w \rangle \nonumber\\
&\le \|\sum_{i\in S}  \epsilon_i  v_i  \|_1 \| 1_E T_n w\|_\infty \nonumber\\
&\le  \|\sum_{i\in S}  \epsilon_i v_i \|_2 \|T\|/\epsilon = |S|^{1/2} \|T\|/\epsilon,
\end{align}
and hence  $|S| \le (\|T\|/\epsilon)^2 (1 - 2\epsilon)^{-2}$.  Therefore  $|\{w_i : \|T w_i\|_1 \ge (1- 2\epsilon)^2\}|$ is at least \break $(1- 2\epsilon)^2 (\epsilon/\|T\|)^2 N(p,n)$.
\qed

\noindent
In the proof of Theorem \ref{theorem} it was not important that the $v_i$ be characters,. One can get similar examples by using the earlier and softer result  \cite{bgn}, \cite{bdgjn}, that for $2<q<\infty$ and every $N$ there are constant modulus vectors $v_1, \dots v_N$ in $L_q^{N^{q/2}}$ so that the $\|\cdot \|_q$ and  $\|\cdot \|_2$ norms are $C_q$ equivalent on the linear span of $(v_i)_{i=1}^N$, and in $L_2^{N^{q/2}}$ every linear combination  of $(v_i)_{i=1}^N$  has its norm  dominated by   the $\ell_2^{N^{q/2}}$ norm of its coefficients. The only change in the proof is that the final equality in (\ref{S}) becomes an inequality.

\noindent
Minor adjustments of the proof of Theorem \ref{theorem}  yield  that other families of closed ideals in $L(L_1)$ are distinct.  For example, for $2<p<\infty$ let $\cJ_p$ be the set of all operators on $L_1$ that factor through an operator of the form $I_{p,1}Si_{1,2}$, where $I_{p,1} : L_p\to L_1$ and $i_{1,2}: \ell_1\to \ell_2$ are the formal identity mappings and $S$ is an operator from $\ell_2$ into $L_p$.  It is easy to check that each $\cJ_p$ is an ideal in $L(L_1)$ and $\cI_p \subset \cJ_p$.  On the other hand, the argument we gave to show that the operator $J_pQ$  (where $Q$ is a quotient mapping from $L_1$ onto $\ell_1$)  is not in  $\overline{\cI_q}$ when $2<p<q<\infty$ also shows that $J_pQ$ is not in $\overline{\cJ_q}$, so  $\overline{\cJ_p}\not=\overline{\cJ_q}$ when $2<p<q<\infty$.  The family of $\cJ_p$'s has some advantages over the family of  $\cI_p$'s: $\cJ_p$ is easy to describe and it is obvious that the family of $\cJ_p$'s are nested.  On the other hand, each $\cI_p$ is generated by a single operator, $J_p$, and with a bit more work we could have constructed the family of $\cI_p$'s to be nested. By the way, we do not know whether the containment $\overline{\cI_p}\subset  \overline{\cJ_p}$ is proper.  In fact, there is a lot of choice in the construction of   $\cI_p$, and we do not know whether different choices produce the same closed ideals.

\noindent
There are known to exist $\aleph_1$ closed ideals in $L(C[0,1])$ since, for example, there are $\aleph_1$ mutually non isomorphic complemented subspaces of $C[0,1]$ each of which is isomorphic to its square \cite{bespel}, but as far as we know, only finitely many small ideals in $L(C[0,1]))$ are known to exist.
To close this section we prove that there are a continuum of small closed ideals in $L(C[0,1])$. Of course, by Miljutin's Theorem \cite{mil}, \cite[Theorem 4.4.8]{albiackalton}, this is equivalent to saying that $L(C(K))$ has a continuum of small closed ideals for some (or all)  uncountable compact metric spaces $K$, so our theorem is proved for the convenient case that $K$ is a Cantor set.

\begin{corollary}\label{C(delta)}
The Banach algebra $L(C(\Delta))$ of  bounded linear operators on the space of continuous functions on the Cantor set $\Delta:=\{-1,1\}^\N$ contains a continuum of small closed ideals.
\end{corollary}
\pf
For $2<p<\infty$ define  $K_p:C(\Delta)\to c_0$ by mapping $f\in C(\Delta)$ to $\int  \chi_n f$, where $\chi_n$ ranges over the characters in the set $V_p$ used in Theorem \ref{theorem}. We can and do identify  $c_0$ with a norm one complemented subspace of $C(\Delta)$.   The enumeration is chosen so that $K_p^* = J_p$, where $J_p$ is the operator from $\ell_1$ into $L_1(\{-1,1\}^\N)$ used in Theorem \ref{theorem} and $L_1(\{-1,1\}^\N)$ is identified with a subspace of $C(\Delta)^*$ in the usual way.  Let ${\cal G}_p$ be the ideal of all operators on 
$C(\Delta)$ that factor through ${ K}_p$.  We proved in Theorem \ref{theorem} that for $2<p<q<\infty$, if $Q$ is a quotient mapping from $L_1(\{-1,1\}^\N)$ onto $\ell_1$, then $J_pQ$ is not in $\overline{\cal I}_q$.  Using the fact that $L_1(\{-1,1\}^\N)$ is a complemented subspace of $C(\Delta)^*$, we deduce  that $K_p $, when considered as an operator in $L(C(\Delta))$ via the identification of $c_0$ with a subspace of $C(\Delta)$,  is not in $\overline{\cal G}_q$.
  Consequently, $\overline{\cal G}_p \not=\overline{\cal G}_q$ when $2<p<q<\infty$. Finally, notice that the ideals  $\overline{\cal G}_p$ are small. Indeed,  $K_p$ factors through a Hilbert space because  $K_p^* = J_p$ does, hence $K_p$  is strictly singular because it ranges in $c_0$.  \endpf

\section{Ideals in $L(\ell_\infty)$}
\label{ellinfinity}
In this section we prove that $L(\ell_\infty)$ has a continuum of closed ideals.  Since $\ell_\infty$ is isomorphic as a Banach space to $L_\infty$,  the Banach algebras $L(\ell_\infty)$ and $L(L_\infty)$ are isomorphic as Banach algebras even though $L(L_1) $ and $L(\ell_1)$ are very different as Banach algebras--$K(\ell_1)$ is the only closed ideal in $L(\ell_1)$ while in Section \ref{ideals} we proved that $L(L_1) $ has a continuum of closed ideals.  The problem in dealing with $L(L_\infty)$ is that for non reflexive spaces $X$, distinct closed ideals in $L(X)$ do not naturally generate distinct closed ideals in $L(X^*)$. For example, $L(L_1)$ has at least two proper closed large ideals; namely, the ideal of operators on $L_1$ that factor through $\ell_1$ and the unique maximal ideal;  while $L(L_\infty)$ has no proper large ideal because the identity on $L_\infty$ factors through every non strictly singular operator on $L_\infty$. However, in Corollary \ref{L1duality} we prove that any family of distinct small closed ideals in $L(L_1)$ naturally ``dualize" to give a corresponding family of small closed ideals in $L(L_\infty)$. This is done in two steps, and only the first step is needed to see that $L(L_\infty)$ has a continuum of closed ideals.  The first step, Proposition \ref{duality}, implies that if $X$ is isomorphic to $X \oplus X$, $X$ is complemented in $X^{**}$,  and $T\in L(X)$ is such that  $TB_X$ is far from $SB_X$ for every $S$ in some  ideal $\cI$ that is contained in the ideal $\cW(X)$ of weakly compact operators on $X$, then   $T^{**}$ is not in the closed ideal in $L(X^{**})$ generated by $\cI^{\#\#}$,
where for a subset $A$ of $L(X)$, $A^{\#}$ is defined to be $\{ S^* \in L(X^*) : S \in A\}$.
It is also convenient to denote by $A^{\dagger}$ the set
$L(X^*)\cdot A^{\#} \cdot L(X^*)$, the collection of all operators in $L(X^*)$ that factor through some operator in $A^\#$. Note  that  $\cI^{\dagger}$ is an ideal in $L(X^*)$ if $\cI$ is an ideal in $L(X)$  and $X \oplus X$ is isomorphic to $X$.  Indeed, then  if $T$, $S$ are in $\cI$, the operator $T\oplus S : X \oplus X \to X \oplus X$ is similar to an operator in the ideal $\cI$, from which it easily follows that $\cI^{\dagger}$  is closed under addition and is an ideal in $L(X^*)$. It is equally easy to verify that if  $\cI$ is an ideal in $L(X)$  and $X\oplus X \sim X$ then $\cI^{\dagger\dagger} = L(X^{**} )  \cdot \cI^{\#\#} \cdot
L(X^{**})$.

\noindent
Given subsets $A$, $B$ of a Banach space $X$, define a (non symmetric) distance from $A$ to $B$ by $d(A,B):= \sup_{x\in A} \inf_{y\in B} \|x-y\|$.

\begin{proposition}\label{duality}
Let $\cI \subset \cW(X)$ be an ideal in $L(X)$, where $X$ is a Banach space that is isomorphic to its square $X\oplus X$. Asssume that there is a projection $P$ from $X^{**}$ onto $X$. Suppose that $T\in L(X)$ and $d(TB_X,SB_X) \ge  \epsilon>0$ for every $S$ in $\cI$. Then $d(T^{**}B_{X^{**}},UB_{X^{**}}) \ge \epsilon/\|P\|$ for every $U$ in $\cI^{\dagger\dagger}$. Consequently, $T^{**}$ is not in the closure of the ideal $\cI^{\dagger\dagger}$ in $L(X^{**})$ and thus $T^*$ is not in the the closure of the ideal $\cI^{\dagger}$ in $L(X^{*})$
\end{proposition}
\pf
Assume, for contradiction, that $d(T^{**}B_{X^{**}},UB_{X^{**}}) < \delta < \epsilon/\|P\|$ for some $U$ in $\cI^{\dagger\dagger}$. Write $U=ES^{**}F$ where $E$, $F$ are norm one operators in $L(X^{**})$ and $S$ is in $\cI$. Thus
\[
TB_X \subset T^{**}B_{X^{**}} \subset ES^{**}B_{X^{**}} + \delta B^o_{X^{**}}
\]
so that
\[
TB_X \subset (PE)S^{**}B_{X^{**}} + \delta \|P\| B^o_{X}.
\]
Since $S$ is weakly compact, the set $SB_{X}$  is norm dense in $S^{**}B_{X^{**}} $.  Therefore
\[
TB_X \subset (PE)SB_X + \delta \|P\| B^o_{X}.
\]
This implies that $d(TB_X,(PE)SB_X) \le \delta \|P\| < \epsilon$, which is a contradiction because $(PE)S$ is in $\cI$.

\noindent The ``consequently" statement is obvious.  \endpf

\begin{theorem}\label{Linfty}
The Banach algebra $L(\ell_\infty)$ has a continuum of distinct closed ideals.
\end{theorem}
\pf
As was mentioned earlier, it is equivalent to prove that $L(L_\infty)$ has a continuum of distinct closed ideals. The space $L_1$ satisfies the hypotheses on the space $X$ in Proposition \ref{duality}.  Let $\cI_p$, $2<p<\infty$, denote the ideals in $L(L_1)$ that were constructed in Theorem \ref{ideals}.  Let $2<p<q$.  In Theorem \ref{ideals} it was proved that $d(V_p, SB_{L_1}) \ge 1/10$ for every $S\in \cI_q$. Now $V_p \subset J_pQB_{L_1}$ and $J_pQ$ is in $\cI_p$, so by Proposition \ref{duality} the operator $(J_pQ)^{**}$ is not in the closure of $\cI_q^{\dagger\dagger}$ and hence $(J_pQ)^*$ is not in the closure of $ \cI_q^{\dagger}$. This proves that the closures of $ \cI_q^{\dagger}$, $2<p<\infty$, are distinct ideals in $L(L_\infty)$. \endpf

\noindent
We turn to the  second step in proving that any family of distinct small closed ideals in $L(L_1)$ gives rise to a family of distinct closed ideals in $L(L_\infty)$.

\begin{proposition} \label{smallL1}
Suppose that $T\in L(L_1)$ and $S\in \cW(L_1)$  and $\epsilon > 0$ are such that \hfill \break $d(TB_{L_1}, SB_{L_1}) < \epsilon$. Then there is $U \in L(L_1)$ such that $\|T - SU\| < \epsilon$.
\end{proposition}
\pf
In the first step we show that there is  $T_1\in L(L_1)$ with $\|T-T_1\| < \epsilon$ and $T_1B_{L_1} \subset \overline{SB_{L_1}}$. In particular,  the operator $T_1$ is weakly compact. The operator $T_1$ is actually of the form $S^{**}V$ for some norm one operator $V:L_1 \to L_1^{**}$ (note that $S^{**}$ ranges in $L_1$ because $S$ is weakly compact).
The proof was more or less  given in the introduction even if the result was not stated there  in this generality, but for convenience we repeat the argument. Write $L_1$ as the closure of the union of a sequence of subspaces $E_n$ with $E_n$ isometrically isomorphic to $\ell_1^n$  and let $P_n$ be a contractive projection from $L_1$ onto $E_n$.  Let $\{e_{n,i}\}_{i=1}^n$ be a basis for $E_n$ that is isometrically equivalent to the unit vector basis of $\ell_1^n$.  Take  $\delta$ so that $ d(TB_{L_1}, SB_{L_1})  <  \delta < \epsilon$ and choose vectors $\{x_{n,i}\}_{i=1}^n$ in $B_{L_1}$ so that $\|  Te_{n,i} - Sx_{n,i}\| < \delta$. Let $A_n$ be $P_n$ followed by the operator extension to $E_n$ of $e_{n,i} \mapsto x_{n.i}$, $1\le i \le n$.  Clearly $\|Tx - SA_nx\|  <\delta$ for every $x\in E_n$.  Considering the   $A_n$ as operators into $L_1^{**}$ with the weak$^*$ topology, we have that some subnet  of $\{A_n\}_n$   converges pointwise weak$^*$ to an operator $V:L_1\to L_1^{**}$. The operator $S^{**}: L_1^{**} \to L_1$ is weak$^*$ to weak continuous and so $\|Tx - S^{**} Vx\| \le \delta $ for every $x\in \cup_{n=1}^\infty E_n$ and hence for every $x\in L_1$.  Since $\|V\| \le 1$ and $S$ is weakly compact, $T_1 := S^{**} V$ maps $B_{L_1}$ into $\overline{SB_{L_1}}$.

\noindent
Before proceeding with the proof, we make a comment. The simple first step gives an operator that factors through $S^{**}$ and  is $\epsilon$-close to $T$.  Perhaps under some reasonable but general conditions this gives an approximation to $T$ that factors through $S$ itself, but we do not see it.  The argument that follows is short but uses quite a bit of the classical structure theory of $L_1$.

\noindent
Since the operator $T_1$ is weakly compact, it is representable \cite[p. 73]{du}. That is, there is a Bochner measurable function $g:(0,1)\to \overline{T_1B_{L_1}}$ so that for every $f\in L_1$ we have $T_1 f = \int_0^1 f(t)g(t) \,dt$ (in this proof we  use $L_1(0,1)$ as our model for $L_1$). Since $g$ is Bochner measurable, it can be approximated arbitrarily closely in the norm of $L_\infty([0,1], L_1)$ by a countably valued Bochner measurable function \cite[p. 42]{du}.  That is, there is a sequence $\{x_n\}_n $ in $g[0,1] \subset \overline{T_1B_{L_1}} \subset \overline{SB_{L_1}}$ and a measurable partition $\{S_n\}_n$ of $[0,1]$ such that the operator $T_2\in L(L_1)$ defined by $T_2f := \sum_{n=1}^\infty (\int_{S_n} f(t) \,dt) x_n$ satisfies the inequality $\|T_1 - T_2\| < \epsilon  - \|T-T_1\|$. The operator $T_2$ has a natural factorization $T_2 = U_1V_1$ through $\ell_1$; indeed, $V_1f := \sum_{n=1}^\infty  (\int_{S_n} f(t) \,dt) e_n$ and $U_1 e_n:= x_n$.  So $\|V_1\| \le 1$ and $U_1 B_{\ell_1} \subset \overline{SB_{L_1}}$. For the final step, take $0 < \tau < \epsilon - \|T-T_1\|- \|T_1 - T_2\|$ and choose $y_n \in B_{L_1}$ so that $\|Sy_n - x_n\| < \tau$. Define an operator $R: \ell_1 \to L_1$ by setting $Re_n = y_n$, so $\|R\| \le 1$.  It is evident that the operator $U:= RV_1$ satisfies the conclusions of  Proposition \ref{smallL1}.
\endpf

\noindent Our next result is an immediate consequence of Proposition \ref{duality} and
Proposition \ref{smallL1}.
\begin{corollary}\label{L1duality}
If $\cI_1$ and  $\cI_2$ are small ideals in $L(L_1)$ that have distinct closures, then the ideals  $\cI_1^\dagger$  and  $\cI_2^\dagger$ in $L(L_\infty)$ also have distinct closures.
\end{corollary}

\bigskip

\bigskip

%\vfill\eject

\noindent William B. Johnson\newline
             Department Mathematics\newline
             Texas A\&M University\newline
             College Station, TX, USA\newline
             E-mail: johnson@math.tamu.edu

\bigskip
\noindent Gilles Pisier\newline
             Department Mathematics\newline
             Texas A\&M University\newline
             College Station, TX, USA\newline
E-mail: pisier@math.tamu.edu

\bigskip
\noindent Gideon Schechtman\newline Department of
Mathematics\newline Weizmann Institute of Science\newline
Rehovot, Israel\newline E-mail: gideon.schechtman@weizmann.ac.il

\end{document}